\UseRawInputEncoding

\documentclass[11pt]{article}
\usepackage{latexsym}
\usepackage{amsmath,amssymb,amstext,amsthm} 
\usepackage{enumerate}
\def\hypergeom#1#2#3#4#5{{}_#1 F_{#2}\left({#3\atop#4}; \ #5\right)}

\def\eqref#1{{\normalfont(\ref{#1})}}

\def\eqref#1{{\normalfont(\ref{#1})}}

\newtheorem{theorem}{Theorem}[section]

\newtheorem{prop}[theorem]{Proposition}

\newtheorem{thm}[theorem]{Theorem}
\newtheorem{cor}[theorem]{Corollary}

\newtheorem{remark}[theorem]{Remark}

\newtheorem{lemma}[theorem]{Lemma}

\newcommand{\A}{{\mathcal A}}

\newcommand{\bbm}{\begin{bmatrix}}
\newcommand{\ebm}{\end{bmatrix}}
\newcommand{\bem}{\begin{pmatrix}}
\newcommand{\eem}{\end{pmatrix}}
\newcommand{\beq}{\begin{equation}}
\newcommand{\beqs}{\begin{equation*}}
\newcommand{\bet}{\begin{table}}
\newcommand{\eeq}{\end{equation}}
\newcommand{\eeqs}{\end{equation*}}
\newcommand{\beqr}{\begin{eqnarray}}

\newcommand{\nc}{\newcommand}
\nc{\arrow}{{\rm arrow\,}}
\nc{\Arrow}{{\rm Arrow\,}}
\nc{\BoDiag}{{\rm B^0Diag\,}}
\nc{\bodiag}{{\rm b^0diag\,}}

\nc{\Mm}{{\mathcal M}^{m} }
\nc{\Mmn}{{\mathcal M}^{mn} }
\nc{\Mnr}{{\mathcal M}_{nr} }
\nc{\Mnmr}{{\mathcal M}_{(n-1)r} }
\nc{\kwqqp}{Q{$^2$}P\,}
\nc{\kwqqps}{Q{$^2$}Ps}

\nc{\notinaho}{(X,S)\in \overline{AHO}(\A)}
\nc{\inaho}{(X,S)\in AHO(\A)}

\newcommand{\bea}{\begin{eqnarray}}%
\newcommand{\eea}{\end{eqnarray}}%
\newcommand{\beas}{\begin{eqnarray*}}%
\newcommand{\eeas}{\end{eqnarray*}}%
{}

\newcommand{\Hnp}[1][]{\,\mathbb{H}_+^{\ifthenelse{\equal{#1}{}}{n}{#1}}}
\newcommand{\Hn}[1][]{\,\mathbb{H}^{\ifthenelse{\equal{#1}{}}{n}{#1}}}
\newcommand{\Dn}[1][]{\,\mathbb{D}^{\ifthenelse{\equal{#1}{}}{n}{#1}}}

\usepackage{tikz}
\usepackage{pgfplots}
\pgfplotsset{compat = newest}
 \usetikzlibrary{arrows.meta,
                backgrounds,
                fit,
                matrix}
                
                \usepackage[symbol]{footmisc}

\begin{document}

\title{The autoregressive filter problem for multivariable degree one symmetric polynomials}
             \author{Jeffrey S. Geronimo\thanks{School of Mathematics, Georgia Institute of Technology, 225 North Ave, Atlanta, GA 30332.},
Hugo J. Woerdeman\thanks{Department of Mathematics, Drexel University, 
3141 Chestnut Street, Philadelphia, PA 19104.  Research supported by 
Simons Foundation grant 355645 and National Science Foundation grant DMS 2000037.
}, and
Chung Y. Wong%
        \thanks{Department of Mathematics, County College of Morris, 214 Center Grove Rd.,
Randolph, NJ 07869.}
}
\date{ }
          \maketitle



\begin{abstract}
The multivariable autoregressive filter problem asks for a polynomial $p(z)=p(z_1, \ldots , z_d)$ without roots in the closed $d$-disk based on prescribed Fourier coefficients of its spectral density function  
 $1/|p(z)|^2$. The conditions derived in this paper for the construction of a degree one symmetric polynomial reveal a major divide between the case of at most two variables vs. the the case of three or more variables. The latter involves multivariable elliptic functions, while the former (due to [J. S. Geronimo and H. J. Woerdeman,  {\it Ann. of Math. (2)}, 160(3):839--906, 2004]) only involve polynomials. The three variable case is treated with more detail, and entails hypergeometric functions. Along the way, we identify a seemingly new relation between $\hypergeom{2}{1}{\frac13,\frac23}{1}{z}$ and 
$\hypergeom{2}{1}{\frac12,\frac12}{1}{\widetilde{z}}$.
\end{abstract}




\section{Introduction}
\label{sec:intro}

The identification problem for wide sense stationary
autoregressive stochastic processes is a classical signal
processing problem. We consider (wide sense)
stationary processes $X_{m}=X_{(m_1, \ldots , m_d)}$ depending on $d$ discrete variables
defined on a fixed probability space $( \Omega , {\mathcal A} , P )$. We
shall assume that the random variables $X_{m}$ are {\em
centered}, i.e., their means $E (X_{m} )$
equal zero. Recall that the space $L^2 ( \Omega , {\mathcal A} , P )$ of
square integrable random variables endowed with the {\it inner
product of centered random variables}
$$ \langle X , Y \rangle := E ( Y^* X ) $$
is a Hilbert space. A sequence $X = ( X_{m} )_{m \in
{\mathbb Z}^d}$ is called a {\em stationary
process} on ${\mathbb Z}^d$ if for
$m,k \in {\mathbb Z}^d$ we have that $$ E ( X_{m}^* X_{k } )
= E( X_{m+p}^* X_{k+p} ) =: R_X(m-k) \
\hbox{\rm for all} \ p \in {\mathbb Z}^d. $$ It is known that the
function $R_X$, termed the {\em covariance
function} of $X$, defines a
positive semidefinite function, that is,
$$ \sum_{i,j=1}^k \alpha_i \overline{\alpha}_j R_X(r_i - r_j) \ge 0$$
for all $k \in {\mathbb N}$, $\alpha_1 , \dots ,\alpha_k \in {\mathbb C},
r_1,\ldots ,r_k \in {\mathbb Z}^d$. Bochner's theorem
\cite{Boch, Boch1} on positive semidefinite functions states that for such
a function $R_X$ there is a positive measure $\mu_X$ defined for
Borel sets on the torus $[0,2 \pi ]^d$ such
that
$$ R_X(r) = \int e^{-i \langle r, u \rangle} d \mu_X (u) $$
for all $d$-tuples of integers $r\in {\mathbb Z}^d$. The measure $\mu_X$ is referred to
as the {\em spectral distribution measure} of the process $X$.

For $n=(n_1, \ldots , n_d) \in{\mathbb N}_0^d$ we let $\underline{n} =\prod_{j=1}^d \{ 0, \ldots , n_j \}$.
A centered stationary stochastic process $X$ is said to be
AR($\underline{n}$) if there exist complex
numbers $a_{k} , k\in \underline{n}
\setminus \{ 0 \}$, such that for every $t$,
\begin{equation}\label{tone} x_{t}+\mathop{\sum\limits_{k\in\underline{n}
}}_{k\ne 0}
 a_{k}x_{t-k}=e_{t}
,\qquad t\in{\mathbb Z}^d, \end{equation} where $\{e_{k}\ ; k
\in{\mathbb Z}^d \}$ is a white noise zero mean process with variance
$\sigma^2$. Here AR stands for
auto-regressive. Let $H$ be the standard half-space in ${\mathbb Z}^d$; that is
$$ H=\{ (k_1, \ldots , k_d) \in {\mathbb Z}^d : {\rm there \ is } \ j\in\{1,\ldots , d \} \ {\rm with} \ k_1=\cdots = k_{j-1} = 0 \ {\rm and} \ k_j >0 \} . $$ The
AR($\underline{n}$) process is said to be {\it
causal} if there is a solution to \eqref{tone}
of the form
\begin{equation}\label{toneb}
x_{t} = \mathop{\sum\limits_{k\in H\cup\{ 0\}}}
\phi_{k}e_{t-k}, \ \ t \in {\mathbb Z}^d , \end{equation} with
$\mathop{\sum_{k\in H\cup\{ 0\}}} \ | \phi_{k}|<\infty$. Causality based on halfspaces and multivariable generalizations of the one variable case go back to the influential papers by Helson and Lowdenslager \cite{HL,HL2}. It
is not difficult to see that the
$AR(\underline{n})$ process $X$ is causal if
and only if the polynomial
$$\tilde{p} (z) = 1+ \mathop{\sum\limits_{k\in\underline{n} }}_{k\ne 0}
\overline a_{k}z^k$$ has no roots in the closed $d$-disk; we call such a polynomial {\it stable}. A causal
AR($\underline{n}$) process is in fact {\it
positive orthant causal}, which by
definition means that there is a solution to \eqref{tone} of the
form
\begin{equation}\label{toneb2}
x_{t} = \mathop{\sum\limits_{k\ge 0}}_{k\ne 0}
\phi_{k}e_{t-k}, \ \ t \in {\mathbb Z}^d , \end{equation}
where $k=(k_1,\ldots , k_d) \ge 0$ means that $k_j\ge 0$ for $j=1,\ldots , d $.

The {\it multivariate autoregressive filter design
problem} is the following. ``Given are covariances
$$c_{k}=E(X_{0}^*X_{k}), \qquad
k\in\underline{n} .$$ What conditions must
the covariances satisfy in order that these are the covariances of
a causal AR($\underline{n})$ process? And in
that case, how does one compute the filter coefficients $a_{k}$,
$k \in \underline{n}\setminus \{ 0 \}$
and $\sigma^2$?"
The papers \cite{Kailath1986},
\cite{Lev-Ari1989}, \cite{Kayran1996} are useful sources for an
explanation how the autoregressive filters are used in signal
processing.

The following characterization for the two variable autoregressive filter design
problem appeared in \cite{GW}.
\begin{theorem}\cite{GW}\label{GW2var}   Let $n,m\in{\mathbb N}$ and $c_{kl}$, $(k,l)\in\{ 0, \ldots , n \}\times\{ 0, \ldots , m \}$, be given complex numbers. There
exists a causal autoregressive process
with the given covariances $c_{kl}$ if and only if there exist
complex numbers $c_{kl}$, $(k,l)\in\{1,\dots,n\} \times
\{-m,\dots, 1\}$, such that
\begin{enumerate} \item the $(n+1)(m+1)\times (n+1)(m+1)$ doubly
indexed Toeplitz matrix $\Gamma=(c_{t-s})_{s,t \in
\{ 0, \ldots , n \}\times\{ 0, \ldots , m \} }$ is positive definite; \item the
matrix $(c_{s-t})_{s\in \{ 1, \ldots , n \} \times \{ 0, \ldots , m \} ,
t \in \{ 0, \ldots , n \}  \times \{ 1, \ldots , m \} }$ has rank equal
to $nm$.
\end{enumerate}
\noindent In this case one finds the vector
$$ \frac{1}{\sigma^2} [ a_{nm}  \cdots a_{n0}  \ \cdots \ a_{0m} \cdots a_{01}
\ 1 ]$$ as the last row of the inverse of $\Gamma$.
\end{theorem}
If we consider the polynomial 
$p(z)= \frac{1}{\sigma}\tilde{p}(z)$, then the Fourier coefficients of $\frac{1}{|{p}|^2}$ coincide exactly with the covariances $c_k$. In other words,
$$ \widehat{\frac{1}{|{p}|^2}} (k) = c_k , k\in\underline{n},  $$
where $\widehat{f} (k)$ denotes the $k$th Fourier coefficient of the function $f$. In the remainder of the paper we will formulate the problems and results in terms of this direct connection. 

In this paper we will focus on the case where the polynomial $p(z)$ is a 
degree one symmetric polynomials in $d$ variables, i.e., 
$$ p(z_1, \ldots , z_d) = p_0 + p_1 (z_1+\cdots + z_d ). $$
In general, a symmetric polynomial is a polynomial where a permutation of the variables does not change the polynomial.
It is easy to see that $p$ is stable if and only if $d|p_1|<|p_0|$. The corresponding autoregressive filter problem is as follows. 

{\bf Problem.} Given $a$ and $b$. Find, if possible, a degree one stable symmetric polynomial in $d$ variables so that 
$$ \widehat{\frac{1}{|p|^2}} (0, 0 , \ldots , 0 ) = a, \widehat{\frac{1}{|p|^2}} (1, 0 , \ldots , 0 ) =b . $$
Clearly, due to the symmetry, we have that 
$$ \widehat{\frac{1}{|p|^2}} (1, 0 , \ldots , 0 ) = \widehat{\frac{1}{|p|^2}} (0,1, \ldots , 0 ) = \cdots = \widehat{\frac{1}{|p|^2}} (0 , \ldots, 0 , 1 ), $$
so that it suffices to just require $\widehat{\frac{1}{|p|^2}} (1, 0 , \ldots , 0 ) =b$.
Notice that $a>0$ will be a necessary condition for the existence of a solution.
If we apply Theorem \ref{GW2var} to this case, we obtain the following.

\begin{theorem} The above problem has a solution in $d=2$ variables if and only if $|b|<a$. In that case, the polynomial $p(z)=p_0 + p_1(z_1+z_2)$ is given via 
$$ \begin{bmatrix} a & \bar{b} & \bar{b} \cr b & a & \frac{|b|^2}{a} \cr b & \frac{|b|^2}{a} & a \end{bmatrix}\begin{bmatrix} |p_0|^2 \cr p_1\bar{p_0}  \cr p_1\bar{p_0} \end{bmatrix} = \begin{bmatrix} 1 \cr 0 \cr 0 \end{bmatrix}. $$ 
\end{theorem}
Indeed, in this case $c_{1,-1}$ is the only unknown in the matrix $\Gamma$, and item 2 in Theorem \ref{GW2var} requires
$$ \begin{bmatrix} c_{1,-1} & c_{0,-1} \cr c_{1,0} & c_{0,0} \end{bmatrix} = \begin{bmatrix} c_{1,-1} & \bar{b} \cr b & a \end{bmatrix}$$ to be of rank 1, which leads to $c_{1,-1}=\frac{|b|^2}{a}$. 

The main result in this paper addresses the case of $d$ variables, which we will state in the next section. Recall that the hypergeometric function $\ _2 F_1$ is defined for $|z|<1$ via the power series
$$ \hypergeom{2}{1}{a,b}{c}{z} =\sum_{n=0}^\infty \frac{(a)_n (b)_n}{(c)_n} \frac{z^n}{n!} .$$ 
Here the Pochhammer function $(q)_n$ is defined by
$$ (q)_n = \begin{cases} 1 , & n=0; \\ q(q+1)\cdots(q+n-1) , & {\rm otherwise} . \end{cases}$$
When we specify the result for $d=3$ variables we obtain the following. 

\begin{theorem}\label{3varcase} The above problem has a solution in $d=3$ variables if and only if $|b|<a$. In that case, one finds the polynomial $p(z)=p_0 + p_1(z_1+z_2+z_3)$ by determining $c\ge 0$ so that 
\begin{equation}\label{hypeq} \frac{a(a+2c)}{a^2+2ac-3|b|^2}=\frac{(a+2c)^2}{(a+2c)^2-3|b|^2}\hypergeom{2}{1}{\frac13,\frac23}{1}{\frac{27|b|^4((a+2c)^2-|b|^2)}{((a+2c)^2-3|b|^2)^3}}, \end{equation} and  
$$ \begin{bmatrix} a & \bar{b} & \bar{b} & \bar{b}\cr b & a &c& c \cr b &c & a & c \cr b & c & c & a \end{bmatrix} $$ is positive definite. Next a solution $p(z)$ is found via the equation
$$ \begin{bmatrix} a & \bar{b} & \bar{b} & \bar{b}\cr b & a &c& c \cr b &c & a & c \cr b & c & c & a \end{bmatrix} \begin{bmatrix} |p_0|^2 \cr p_1\bar{p_0}  \cr p_1\bar{p_0} \cr p_1\bar{p_0} \end{bmatrix} = \begin{bmatrix} 1 \cr 0 \cr 0 \cr 0 \end{bmatrix}. $$\end{theorem}

As one can see there is a significant difference between two and three variables. In two variables the unknown in the matrix is easily found by setting $c=\frac{|b|^2}{a}$, while in three variables one needs to solve the highly nontrivial equation \eqref{hypeq} to find the unknown $c$ in the matrix.  The number $c$ plays the role of $$ c= \widehat{\frac{1}{|p|^2}} (1,-1,0,\ldots ,0) = \widehat{\frac{1}{|p|^2}} (-1,1,0,\ldots ,0) =\cdots = \widehat{\frac{1}{|p|^2}} (0,\ldots ,0,1,-1), $$  where again we used the symmetry of the polynomial. We will see that $c$ is required to be nonnegative (see Proposition \ref{prop22}).

The paper is organized as follows. In Section 2 we present our main result giving necessary and sufficient condition for the existence of an autoregressive filter with a stable symmetric degree one polynomial in $d$ variables, as well as a method how to find the polynomial. In Section 3 we further specify the results for the case of three variables and present a new relation between $\hypergeom{2}{1}{\frac13,\frac23}{1}{z}$ and 
$\hypergeom{2}{1}{\frac12,\frac12}{1}{\widetilde{z}}$. Finally, in Section \ref{sec4} we explore finding formulas for other Fourier coefficients in the three variable case. 

%
%
%
%
%
%
%
%
%
%
%
%

\section{The main result}

We will begin by determining some of the Fourier coefficients of $\frac{1}{|p(z)|^2}$, where $p(z)=p_0 + p_1 (z_1+\cdots + z_d)$, $z=(z_1, \ldots z_d)$. It will be convenient to do a simple scaling and assume that $p_0=1$. Next we will  write $p_1 = -s$. We will use the notation
$$ {\mathbb D} = \{ z \in {\mathbb C}: |z|<1 \} ,  {\mathbb T} = \{ z \in {\mathbb C}: |z|=1 \}, 
\overline{{\mathbb D}} = {\mathbb D} \cup {\mathbb T}, {\mathbb N}_0 = \{ 0,1,2,\ldots \} . $$ 

\begin{lemma} The polynomial $ p(z) = 1 -s(z_1+\cdots + z_d) $ is stable if and only if $|s|<\frac{1}{d}$. \end{lemma}

{\it Proof.} Let $|s|<\frac{1}{d}$ and $(z_1, \ldots z_d) \in\overline{\mathbb D}^d$. Then we have that $|s(z_1+\cdots + z_d)|<1$, and thus $p(z)\neq 0$. This gives that $p(z)$ is stable.

When $|s|\ge\frac{1}{d}$, then $z_1=\cdots =z_d = \frac{1}{sd}$ yields a root of $p(z)$ inside $\overline{\mathbb D}^d$. Thus $p(z)$ is not stable. \hfill $\square$
  
For $q\in{\mathbb Z}$ we let $q^+=\max\{ 0, q \}$ and $q^-=\max \{ 0, -q \}$.    
\begin{prop}\label{prop22} Let $ p(z) = 1 -s(z_1+\cdots + z_d)$, $|s|<\frac{1}{d}$. Then for $k=(k_1,\ldots , k_d) \in{\mathbb Z}^d$, 
$$ \widehat{\frac{1}{|p|^2}} (k) = \sum_{n=0}^\infty \sum_{\sum n_i=n} 
\begin{pmatrix} n+k_1^+ +\cdots + k_d^+ \cr n_1+k_1^+ , \ldots , n_d+k_d^+ \end{pmatrix} 
\begin{pmatrix} n+k_1^- +\cdots + k_d^- \cr n_1+k_1^- , \ldots , n_d+k_d^- \end{pmatrix} |s|^{2n}s^{\sum_j k_j^+}\bar{s}^{\sum_j k_j^-} . $$ Here $n_1, \ldots , n_d \ge 0$ range over all nonnegative numbers that sum up to $n$. In particular, 
\begin{equation}\label{ac} \widehat{\frac{1}{|p|^2}} (0, \ldots , 0 ) >0, \widehat{\frac{1}{|p|^2}} (1,-1,0,\ldots , 0) \ge 0. \end{equation}
\end{prop}

{\it Proof.} For $(z_1,\ldots , z_d) \in{\mathbb T}^d$ we have 
$$ \frac{1}{p(z)} = \sum_{n=0}^\infty s^n(z_1+\cdots + z_d)^n = \sum_{n=0}^\infty \sum_{\sum n_i=n} \begin{pmatrix} n \cr n_1, \ldots , n_d  \end{pmatrix} s^{n} z_1^{n_1}\cdots z_d^{n_d}, $$ $$
\frac{1}{\overline{p(z)}}=  \sum_{n=0}^\infty \bar{s}^n(z_1^{-1}+\cdots + z_d^{-1})^n = \sum_{n=0}^\infty \sum_{\sum n_i=n}  \begin{pmatrix} n \cr n_1, \ldots , n_d \end{pmatrix} \bar{s}^{n}  z_1^{-n_1}\cdots z_d^{-n_d}. $$
Multiplying the two and extracting the coefficient of $z_1^{k_1}\cdots z_d^{k_d}$ gives the stated formula for $ \widehat{\frac{1}{|p|^2}} (k)$.

Finally, when $k=(0,\ldots , 0)$ the number $s$ only appears in $|s|^{2n}$ which is always $\ge 0$, and $>0$ when $n=0$, and when $k=(1,-1,0,\ldots , 0)$ the number $s$ only appears in  $|s|^{2n+2}$ which is always $\ge 0$. Clearly, all the multinomial coefficients are nonnegative, and thus \eqref{ac} follows. \hfill $\square$

\begin{prop}\label{prop23} Let $ p(z) = 1 -s(z_1+\cdots + z_d) $, $|s|<\frac{1}{d}$. Put $$a=\widehat{\frac{1}{|p|^2}} (0, 0 , \ldots , 0 ) , b= \widehat{\frac{1}{|p|^2}} (1, 0 , \ldots , 0 ) ,
c= \widehat{\frac{1}{|p|^2}} (1, -1 ,0, \ldots , 0 ) . $$
Then $a>0, c\ge 0$, and the matrix \begin{equation}\label{Amatrix3} A=  \begin{bmatrix} a & \bar{b} & \bar{b} & \cdots & \bar{b} \cr b & a & c & \cdots & c \cr
b & c & a & \cdots & c \cr \vdots & \vdots & & \ddots & \vdots \cr 
b & c & c & \cdots & a \end{bmatrix}, \end{equation}
is positive definite. 
Furthermore 
\begin{equation}\label{Aeq}  
A \begin{bmatrix} 1\cr -s \cr \vdots \cr -s \end{bmatrix} =  \begin{bmatrix} 1\cr 0 \cr \vdots \cr 0 \end{bmatrix}. \end{equation}
\end{prop} 

{\it Proof.} 
Let $$\frac{1}{|p(z)|^2} = \sum_{k\in{\mathbb Z}^d} c_k z^k $$ denote its Fourier series. Thus $\widehat{\frac{1}{|p(z)|^2} } = c_k$, $k\in{\mathbb Z}^d$. Since $\frac{1}{|p(z)|^2}$ is positive, the multiplication operator on $L_2({\mathbb T}^d)$ with symbol $\frac{1}{|p(z)|^2}$ is positive definite. Its matrix representation with respect to the standard monomial basis is $(c_{k-\ell})_{k,\ell\in{\mathbb Z}^d}$. Consequently, any principal submatrix $(c_{k-\ell})_{k,\ell\in \Lambda}$, $\Lambda \subseteq {\mathbb Z}^d$, is positive definite. If we let $\Lambda = \{ 0,e_1, \ldots , e_d \}$, where $e_j$ is the $j$th standard basis vector of ${\mathbb C}^d$, we obtain
\begin{equation}\label{Amatrix4} (c_{k-\ell})_{k,\ell\in \Lambda} =  \begin{bmatrix} a & \bar{b} & \bar{b} & \cdots & \bar{b} \cr b & a & c & \cdots & c \cr
b & c & a & \cdots & c \cr \vdots & \vdots & & \ddots & \vdots \cr 
b & c & c & \cdots & a \end{bmatrix}, \end{equation}
where $$a=\widehat{\frac{1}{|p|^2}} (0, 0 , \ldots , 0 ) , b= \widehat{\frac{1}{|p|^2}} (1, 0 , \ldots , 0 ) ,
c= \widehat{\frac{1}{|p|^2}} (1, -1 ,0, \ldots , 0 ) . $$
Thus \eqref{Amatrix4} is positive definite.  

Next,  we have that 
$$\frac{1}{|p(z_1, \ldots , z_d )|^2} p(z_1, \ldots , z_d) = \frac{1}{ p(\frac{1}{z_1}, \ldots , \frac{1}{z_d} )}= \sum_{k\in{\mathbb N}_0^d} \phi_k z^{-k} , z\in{\mathbb T}^d, $$
where $\phi_0=1$. Comparing the coefficients of $1, z_1, \ldots , z_d$ on both sides we get the equality \eqref{Aeq}. 
\hfill $\square$

\begin{prop}\label{prop24} Let $p_s(z)=1 -s(z_1+\cdots + z_d)$, $|s|<\frac{1}{d}$. Put $$a(s)=\widehat{\frac{1}{|p_s|^2}} (0, 0 , \ldots , 0 ) , b(s)= \widehat{\frac{1}{|p_s|^2}} (1, 0 , \ldots , 0 ) . $$
Then $a(s)$ is a function of $|s|$ and strictly increasing function for $|s|\in [0,\frac{1}{d})$, and $$\{ a(s) : |s|\in [0,\frac{1}{d}) \} = [1,\gamma_d) ,$$ where 
\begin{equation}\label{gamma} \gamma_d = \sum_{n=0}^\infty \sum_{\sum n_i=n} \begin{pmatrix} n \cr n_1, \ldots , n_d \end{pmatrix}^2 d^{-2n} . \end{equation}
We have $\gamma_1=\gamma_2=\gamma_3=\infty$ and $\gamma_d <\infty$ for $d\ge 4$.
Finally, $$ \{ \frac{|b(s)|}{a(s)} :  |s|\in [0,\frac{1}{d}) \} = [0,1-\frac{1}{\gamma_d}), $$ 
where $\frac{1}{\infty}=0$.
\end{prop} 

{\it Proof.} By the established asymptotic that was first ascertained in \cite{RR} and later generalized by \cite[Theorem 4]{RS} and \cite[Theorem 5.1]{BW}, we have 
$$\sum_{\sum n_i=n} \begin{pmatrix} n \cr n_1, \ldots , n_d \end{pmatrix}^2 d^{-2n}
\approx d^{d/2}(4\pi n )^{(1-d)/2} = \Theta (n^{(1-d)/2}) \ \ {\rm as }\ n\to \infty . $$ Thus $\gamma_d=\infty$ for $d\le 3$, and $\gamma_d <\infty$ for $d>3$ follows.
By  Proposition \ref{prop22} we have that
$$ a(s) = \sum_{n=0}^\infty \sum_{\sum n_i =n} \begin{pmatrix} n \cr n_1 , \ldots , n_d \end{pmatrix}^2 |s|^{2n}, $$
thus $a(s)$ is a continuous function and is increasing as $|s|$ increases. 
Further, $a(0)=1$ and $\lim_{|s|\to \frac{1}{d}-} a(s)=\gamma_d$, yielding that the range of $a(s)$ is the interval $[1,\gamma_d )$. 
Similarly, 
$$ |b(s)| =  \sum_{n=0}^\infty \sum_{\sum n_i =n} \begin{pmatrix} n+1 \cr n_1+1 , n_2, \ldots , n_d \end{pmatrix} \begin{pmatrix} n \cr n_1 , \ldots , n_d \end{pmatrix} |s|^{2n+1} $$ is a continuous function and is increasing as $|s|$ increases. Also, note that $b(0)=0$. By \eqref{Aeq} we have that
$$ a(s) - ds\overline{b(s)} =1, $$
and thus 
$$ \frac{|b(s)|}{a(s)}  = \frac{1}{d|s|}\frac{a(s)-1}{a(s)}= \frac{1}{d|s|} (1-\frac{1}{a(s)}) . $$ Since 
$\frac{|b(0)|}{a(0)}=0$ and $\lim_{|s|\to \frac{1}{d}-}  \frac{|b(s)|}{a(s)} = 1-\frac{1}{\gamma_d}, $
the last statement follows. \hfill $\square$

The main result is the following.

\begin{theorem}\label{main} Let $d\ge 3$ and define $\gamma_d$ via \eqref{gamma}. Given are $a>0$ and $b\in{\mathbb C}$ Then there exists a stable degree one symmetric polynomial $p(z_1, \ldots , z_d)$ so that $$ \widehat{\frac{1}{|p|^2}} (0, 0 , \ldots , 0 ) = a, \widehat{\frac{1}{|p|^2}} (1, 0 , \ldots , 0 ) =b , $$
if and only if $|b|<(1-\frac{1}{\gamma_d})a$. In that case, the polynomial $p(z)$ may be found by finding $c\ge 0$ so that
\begin{equation}\label{maineq} \frac{a(a+(d-1)c)}{a^2+(d-1)ac+d|b|^2} =\frac{1}{(2\pi)^{d-2}}  \int_{[0,2\pi]^{d-2}} \frac{1}{\sqrt{g(t_3,\ldots , t_d)}} dt_3\cdots dt_d , \end{equation}
where 
$$ g(t_3,\ldots , t_d) = \left( 1-\frac{2|b|}{a+(d-1)c} \sum_{3\le j \le d} \cos t_j + \frac{|b|^2}{(a+(d-1)c)^2}
\sum_{3\le j,k \le d} \cos (t_j-t_k) \right) \times $$ $$ \ \ \ \ \ \ \ \ \ \ \ \ \ \left(1-\frac{2|b|}{a+(d-1)c} \sum_{3\le j \le d} \cos t_j + \frac{|b|^2}{(a+(d-1)c)^2} \left(
-4+\sum_{3\le j,k \le d} \cos (t_j-t_k) \right) \right) , $$
and the matrix 
$$ \begin{bmatrix} a & \bar{b}  & \bar{b}  & \cdots & \bar{b}  \cr b & a & c & \cdots & c \cr
b & c & a & \cdots & c \cr \vdots & \vdots & & \ddots & \vdots \cr 
b & c & c & \cdots & a \end{bmatrix}$$
is positive definite.
Subsequently, $p(z)=p_0+p_1(z_1+\cdots + z_d)$ is found via the equation
$$ \begin{bmatrix} a & \bar{b}  & \bar{b}  & \cdots & \bar{b}  \cr b & a & c & \cdots & c \cr
b & c & a & \cdots & c \cr \vdots & \vdots & & \ddots & \vdots \cr 
b & c & c & \cdots & a \end{bmatrix} \begin{bmatrix} |p_0|^2 \cr p_1\bar{p_0} \cr \vdots \cr p_1\bar{p_0}  \end{bmatrix} = \begin{bmatrix} 1 \cr 0 \cr \vdots \cr 0 \end{bmatrix}. $$
\end{theorem}

\begin{remark} When we put  $s=\frac{b}{a^2+(d-1)ac-d|b|^2}$, the right hand side of \eqref{maineq} may be rewritten as
 $$  \frac{1}{(2\pi i)^{d-2}} \int_{{\mathbb T}^{d-2}}  \frac{1}{|1-s(z_3+\cdots +z_d)|}  \frac{1}{\sqrt{|1-s(z_3+\cdots +z_d)|^2-4|s|^2}} \ \frac{dz_3}{z_3} \cdots \frac{dz_d}{z_d}. $$
 \end{remark}

In determining the Fourier coefficients of $\frac{1}{|p(z)|^2}$, where 
$$ p(z) = 1 -s(z_1+\cdots + z_d), \ |s|<\frac{1}{d},  $$
we let $w=z_3+\cdots + z_d$, which we will treat as a parameter, and write 
$$ p(z)= p(z_1,z_2,w) = p_{0}(w) -s(z_1 + z_2), $$
where $p_{0}(w) = 1-sw$. 
We write $f(z)=\frac{1}{|p(z)|^2}$ in Fourier series with $w$ as a parameter 
$$ f(z) = \sum_{k,l\in {\mathbb Z}} c_{kl}(w)z_1^kz_2^l . $$

\begin{prop}\label{KWformula} Let $p(z)= p(z_1,z_2,w) = p_{0}(w) -s(z_1 + z_2)$, $p_{0}(w) = 1-sw$, $|s|<\frac{1}{d}$, and write $f(z)=\frac{1}{|p(z)|^2}$ in Fourier series as
$$ f(z) = \sum_{k,l\in {\mathbb Z}} c_{kl}(w)z_1^kz_2^l . $$ Then \begin{equation}\label{3by3inverse}\begin{bmatrix} c_{00}(w) & c_{0,-1}(w) & c_{-1,0}(w) 
\cr c_{01}(w) & c_{00}(w) & c_{-1,1}(w) \cr
c_{10}(w) & c_{1,-1}(w) & c_{00}(w) \end{bmatrix}^{-1} =  \end{equation} $$\begin{bmatrix} |1-sw|^2 &  -\bar{s} (1-sw)  &  -\bar{s} (1-sw)  \cr   -s (1-\bar{s}\bar{w})& \frac12 ( |1-sw|^2+\sqrt{ |1-sw|^4 -4|s|^2|1-sw|^2}\ ) & 0 \cr -s (1-\bar{s}\bar{w}) & 0 &  \frac12 ( |1-sw|^2+\sqrt{ |1-sw|^4 -4|s|^2|1-sw|^2}\ ) \end{bmatrix}  $$
and
$$\begin{bmatrix} c_{00}(w) & c_{0,-1}(w) & c_{-1,0}(w) &c_{-1,-1}(w)
\cr c_{01}(w) & c_{00}(w) & c_{-1,1}(w) & c_{-1,0}(w)\cr
c_{10}(w) & c_{1,-1}(w) & c_{00}(w) & c_{0,-1}(w) \cr 
c_{11}(w) & c_{10}(w) & c_{01}(w) &c_{00}(w)
\end{bmatrix}^{-1} =  $$ $$\tiny \begin{bmatrix} |1-sw|^2 &  -\bar{s} (1-sw)  &  -\bar{s} (1-sw) & 0 \cr   -s (1-\bar{s}\bar{w})& s^2+ \frac12 ( |1-sw|^2+\sqrt{ |1-sw|^4 -4|s|^2|1-sw|^2}\ ) & s^2 &  -\bar{s} (1-sw)\cr -s (1-\bar{s}\bar{w}) & s^2 &  s^2+\frac12 ( |1-sw|^2+\sqrt{ |1-sw|^4 -4|s|^2|1-sw|^2}\ ) &  -\bar{s} (1-sw)\cr 0 &  -s (1-\bar{s}\bar{w}) &  -s (1-\bar{s}\bar{w}) & |1-sw|^2 \end{bmatrix}. $$ 
\end{prop}
 
{\it Proof.} The first inverse follows from \cite[Theorem 1.1]{KW}. With $p(z_1,z_2) = p_{00}+p_{01}z_2+p_{10}z_1+p_{11}z_1z_2$ and using the notation from  \cite[Theorem 1.1]{KW} we have
$$ A=\begin{bmatrix} p_{00} & 0 & 0 \cr p_{01} & p_{00} & 0 \cr p_{10} & 0 & p_{00} \end{bmatrix} , B=\begin{bmatrix} p_{11} & p_{10} & p_{01} \cr 0 & p_{11} & 0 \cr 0 & 0 & p_{11} \end{bmatrix}, $$
$$ C_1= \begin{bmatrix} 0 & 0 & p_{10}\overline{p_{00}}-\overline{p_{01}}p_{11} \cr 0 & 0 & 0 \cr \vdots & \vdots & \vdots \end{bmatrix}, C_2= \begin{bmatrix} 0 &  p_{01}\overline{p_{00}}-\overline{p_{10}}p_{11} & 0 \cr 0 & 0 & 0 \cr \vdots & \vdots & \vdots \end{bmatrix}, $$
$$ D_1 = \begin{bmatrix} |p_{00}|^2 +|p_{10}|^2- |p_{01}|^2 & p_{00}\overline{p_{10}} & 0 & 0 & \cdots \cr
p_{10} \overline{p_{00}} & |p_{00}|^2 +|p_{10}|^2- |p_{01}|^2 & p_{00}\overline{p_{10}} & 0 & \cdots  \cr 
0 &p_{10} \overline{p_{00}} & |p_{00}|^2 +|p_{10}|^2- |p_{01}|^2 & \ddots & \ddots 
\cr \vdots &  \ddots & \ddots  & \ddots & \ddots \end{bmatrix} , $$
$$ D_2 = \begin{bmatrix} |p_{00}|^2 +|p_{01}|^2- |p_{10}|^2 & p_{00}\overline{p_{01}} & 0 & 0 & \cdots \cr
p_{01} \overline{p_{00}} & |p_{00}|^2 +|p_{01}|^2- |p_{10}|^2 & p_{00}\overline{p_{01}} & 0 & \cdots  \cr 
0 &p_{01} \overline{p_{00}} & |p_{00}|^2 +|p_{01}|^2- |p_{10}|^2 & \ddots & \ddots 
\cr \vdots &  \ddots & \ddots  & \ddots & \ddots \end{bmatrix} .$$
To invert $D_1$ we write $D_1= K_1 K_1^*$, where $K_1$ is an upper triangular bidiagonal Toeplitz operator with $\alpha$ on the main diagonal and $\beta$ on the superdiagonal, where $\alpha>0$ and $\beta$ are so that 
$$ \alpha^2 + |\beta|^2= ( |p_{00}|^2 +|p_{10}|^2- |p_{01}|)^2 + | p_{00}\overline{p_{10}}|^2 , ab =  (p_{00}\overline{p_{10}})( |p_{00}|^2 +|p_{10}|^2- |p_{01}|^2). $$ Similarly for $D_2$. Now we use the formula
$$ \begin{bmatrix} c_{00}(w) & c_{0,-1}(w) & c_{-1,0}(w) 
\cr c_{01}(w) & c_{00}(w) & c_{-1,1}(w) \cr
c_{10}(w) & c_{1,-1}(w) & c_{00}(w) \end{bmatrix}^{-1} =AA^*- B^*B - C_1^* D_1^{-1} C_1 - C_2^* D_2^{-1} C_2
$$ to obtain \eqref{3by3inverse}.

For the second inverse, we use that the (4,1) entry in the inverse is 0 as  $p(z)$ does not have a $p_{11}z_1z_2$ term. It now follows from the inverse block matrix formula
\begin{equation}\label{3} \begin{bmatrix} P & H_1 & H_3 \cr H_1^* & Q & H_2 \cr H_3^* & H_2^* & R \end{bmatrix}^{-1}  =
\begin{bmatrix} \begin{bmatrix} P & H_1 \cr H_1^* & Q \end{bmatrix}^{-1} &  \begin{matrix} 0 \cr 0 \end{matrix}  \cr  \begin{matrix} 0 & & 0 \end{matrix} & 0  \end{bmatrix} + \begin{bmatrix} 0 &  \begin{matrix} 0 & & 0 \end{matrix} \cr  \begin{matrix} 0 \cr 0 \end{matrix} & \begin{bmatrix} Q & H_2  \cr H_2^* & R \end{bmatrix}^{-1} \end{bmatrix} - \begin{bmatrix} 0 & 0 & 0 \cr 0 & Q^{-1}  & 0  \cr 0 & 0 & 0   \end{bmatrix},\end{equation} which holds if there is a zero in the (3,1) block of the inverse. 
\hfill $\square$

\begin{prop}\label{prop27} For $p(z)=1 -s(z_1+\cdots + z_d)$, $|s|<\frac{1}{d}$, we have 
$$\widehat{\frac{1}{|p|^2}} (0,  \ldots , 0 )  = \frac{1}{(2\pi i)^{d-2}} \int_{{\mathbb T}^{d-2}}  \frac{1}{|1-s(z_3+\cdots +z_d)|}  \frac{1}{\sqrt{|1-s(z_3+\cdots +z_d)|^2-4|s|^2}} \ \frac{dz_3}{z_3} \cdots \frac{dz_d}{z_d}. $$
\end{prop}

{\it Proof.} In general we have that
$$ \begin{bmatrix} x & \bar{y} & \bar{y} \cr y & v & 0 \cr y & 0 & v \end{bmatrix}^{-1} = \frac{1}{xv-2|y|^2}\begin{bmatrix} v &-\bar{y} & -\bar{y} \cr -y & x-\frac{|y|^2}{v} & \frac{|y|^2}{v} \cr -y & \frac{|y|^2}{v} & x-\frac{|y|^2}{v} \end{bmatrix} . $$ Combining this with Proposition \ref{KWformula} we find $$c_{00}(w) = \frac{v}{xv-2|y|^2},$$
where $$v= \frac12 ( |1-sw|^2+\sqrt{ |1-sw|^4 -4|s|^2|1-sw|^2}\ ) , x=|1-sw|^2, y=-s (1-\bar{s}\bar{w}).$$
We have
$$ xv-2|y|^2= \frac12\left( |1-sw|^4 -4|s|^2|1-sw|^2+ |1-sw|^2 \sqrt{ |1-sw|^4 -4|s|^2|1-sw|^2}\right) = $$
$$  v \sqrt{ |1-sw|^4 -4|s|^2|1-sw|^2} = v |1-sw| \sqrt{ |1-sw|^2 -4|s|^2}.$$
Thus 
$$ c_{00}(w) = \frac{1}{ |1-sw| \sqrt{ |1-sw|^2 -4|s|^2}} = 
\frac{1}{|1-s(z_3+\cdots +z_d)|}  \frac{1}{\sqrt{|1-s(z_3+\cdots +z_d)|^2-4|s|^2}}. $$ To find the 
$0$th Fourier coefficient of $\frac{1}{|p(z)|^2}$ we need to compute 
$$ \frac{1}{(2\pi i)^{d-2}}\int_{{\mathbb T}^{d-2}} c_{00}(z_3+\cdots + z_d) \frac{dz_3}{z_3} \cdots \frac{dz_d}{z_d} , $$ which yields the stated formula. \hfill $\square$

It is easy to check the following lemma.
\begin{lemma}\label{invlemma} Suppose that the $(d+1)\times (d+1)$ matrix 
\begin{equation}\label{Amatrix2} A = \begin{bmatrix} a & \bar{b} & \bar{b} & \cdots & \bar{b} \cr b & a & c & \cdots & c \cr
b & c & a & \cdots & c \cr \vdots & \vdots & & \ddots & \vdots \cr 
b & c & c & \cdots & a \end{bmatrix}\end{equation}
is invertible. Then the first column of the inverse equals 
$$ \frac{1}{a^2+(d-1)ac - d|b|^2} \begin{bmatrix} a+(d-1)c \cr -b \cr \vdots \cr -b \end{bmatrix} . $$
\end{lemma}
{\it Proof.} Simply multiply $A$ by the vector to obtain the first standard basis vector. \hfill $\square$

%
%
%
%

\bigskip

{\it Proof of Theorem \ref{main}.} By the last statement in Proposition \ref{prop24} we see that $\frac{|b|}{a} \in [0, 1-\frac{1}{\gamma_d})$ is necessary and sufficient. 

Next, the polynomial $p(z)$ after normalization so that $p(0)=1$ will satisfy \eqref{Aeq}. Starting with $A$ as in \eqref{Amatrix2} we can, by Lemma \ref{invlemma}, rescale the matrix as $\frac{a+(d-1)c}{a^2+(d-1)ac-d|b|^2} A$ so that the (1,1) entry of its inverse is 1, which corresponds to the situation where $p(0)=1$. Then, again using Lemma \ref{invlemma}, we find that $s=-\frac{dp}{dz_1}|_{z=0} $ corresponds to the value $s=\frac{b}{a^2+(d-1)ac-d|b|^2}$. Using this value for $s$ as well as
$\widehat{\frac{1}{|p|^2}} (0,\ldots, 0 ) = a\frac{a+(d-1)c}{a^2+(d-1)ac-d|b|^2}$, we find that Proposition \ref{prop27} yields equality \eqref{maineq}. \hfill $\square$

%
%
%
%
%

\section{The three variable case}

In this section we provide further details when $d=3$. 
To be consistent with earlier results in \cite{GWW} and \cite{Wong}, we consider the polynomial $$p(z_1,z_2,z_3)= 1 - \frac{z_1+z_2+z_3}{r} , r>3 . $$
Comparing this with the previous section, we make the conversion $s=\frac{1}{r}$ and require $s>0$. This is not a significant restriction as a phase appearing in $s$ can always be absorbed in the variables via $(z_1,z_2,z_3) \to e^{i\theta} (z_1,z_2,z_3)$.

We will use the complete elliptic integral of the first kind, which is
$$K(m)= \int_0^{\frac{\pi}{2}} \frac{1}{\sqrt{1-m\sin^2(t)}}dt = \int_0^1 \frac{1}{\sqrt{1-t^2}\sqrt{1-mt^2}}dt  = \frac{\pi}{2} \hypergeom{2}{1}{\frac12,\frac12}{1}{m}.$$ 
\bigskip

\begin{theorem}\label{thmc000} 
Let $p(z_1,z_2,z_3) = 1-\frac{z_1+z_2+z_3}{r}$, $r>3$, and 
$ f(z)=\frac{1}{|p(z)|^2}$, $z=(z_1,z_2,z_3)$. Write
$$  f(z) =\sum_{k,l,m\in {\mathbb Z}}c_{klm}z_1^k
z_2^l z_3^m , \ \ (z_1,z_2,z_3) \in {\mathbb T}^3  . $$
Then
\begin{equation}\label{c000c1}
 c_{000} =  \frac{r^2}{2\pi}\int_0^{2\pi}  \frac{1}{\sqrt{r^2+1-2r\cos t} \sqrt{ r^2-3-2r\cos t}}  dt= \end{equation}
$$ \ \ \ \ \ \frac{2r^2}{\pi (r-1)^{\frac32}(r+3)^{\frac12}} K(\frac{16r}{(r-1)^3(r+3)}) =\frac{r^2}{(r-1)^{\frac32}(r+3)^{\frac12}}  \hypergeom{2}{1}{\frac12,\frac12}{1}{\frac{16r}{(r-1)^3(r+3)}} .$$
 \end{theorem}

{\it Proof of Theorem \ref{thmc000}.}
From Proposition \ref{prop27} with $s=\frac{1}{r}$ and $z_3=e^{it}$ we obtain
$$c_{000} = \frac{1}{2\pi} \int_0^{2\pi}  \frac{1}{|1-\frac{e^{it}}{r}|}  \frac{1}{\sqrt{|1-\frac{e^{it}}{r}|^2-\frac{4}{r^2}}} dt. $$
Using that $|1-\frac{e^{it}}{r}|^2=(1-\frac{\cos t}{r})^2+\frac{\sin^2 t }{r^2} = \frac{1}{r^2}(r^2 -2r\cos t + 1)$, formula \eqref{c000c1} follows.

Next, use $\cos t = 2\cos^2 \frac{t}{2}-1 = 2\sin^2(\frac{\pi}{2}-\frac{t}{2} )-1$, do a change of variable $t\to \frac{\pi}{2}-\frac{t}{2}$, use the symmetry of the integrand, and \eqref{c000c1} becomes
\begin{equation}\label{GRproofa}  \frac{2r^2}{\pi}  \int_0^{\frac{\pi}{2}}  \frac{1}{\sqrt{(r+1)^2-4r\sin^2 t} \sqrt{ (r+3)(r-1)-4r\sin^2 t}}  dt 
   .  \end{equation}
   Now we let $p^2=\frac{4r}{(r+1)^2}$ and $q^2 = \frac{4r}{(r+3)(r-1)}$, and use 
the first formula in Section 2.616 of \cite{GR}, which is the equality\footnote{due to a change of variables $\sin \alpha = \frac{\sqrt{1-p^2}\sin x}{\sqrt{1-p^2\sin^2x}}$.} 
  $$ \int_0^{\frac{\pi}{2}} \frac{dx}{\sqrt{(1-p^2\sin^2 x) (1-q^2 \sin^2 x)}} = \frac{1}{\sqrt{1-p^2}} \int_0^{\frac{\pi}{2}} \frac{d\alpha}{\sqrt{1-\frac{q^2-p^2}{1-p^2}\sin^2 \alpha}} . $$
  This transforms \eqref{GRproofa} into
  $$  \frac{2r^2}{\pi } \frac{1}{\sqrt{(r-1)^3(r+3)}}\int_0^{\frac{\pi}{2}}  \frac{1}{\sqrt{1-\frac{16r}{(r-1)^3(r+3)}\sin^2 t} }  dt = \frac{2r^2}{\pi (r-1)^{\frac32}(r+3)^{\frac12}} K(\frac{16r}{(r-1)^3(r+3)}) . $$
\hfill $\square$

The following result is inspired by a generating function entry by Paul D. Hanna \cite{Hanna} regarding sequence A002893 on the On-Line Encyclopedia of Integer Sequences (oeis.org). Hanna arrived at this entry as a variation of the generating function for the triangle of cubed binomial coefficients (sequence A181543 on oeis.org) and numerically verified it for hundreds of terms \cite{Hanna2}.

\begin{thm}\label{c000formula}  
Using the same notation as in Theorem \ref{thmc000},
we have
$$  c_{000} = \frac{r^2}{r^2-3} \hypergeom{2}{1}{\frac13 , \frac23}{1}{\frac{27(r^2-1)}{(r^2-3)^3}}, r>3. $$

\end{thm}

{\it Proof.} By Proposition \ref{prop22}, we have $c_{000}(r) = \sum_{n=0}^\infty \sum_{n_1+n_2+n_3=n} \begin{pmatrix} n \cr n_1,n_2,n_3 \end{pmatrix}^2 r^{-2n} , r>3. $
Letting $x=r^{-2}$, and 
$$g(x) = \sum_{n=0}^\infty \sum_{n_1+n_2+n_3=n} \begin{pmatrix} n \cr n_1,n_2,n_3 \end{pmatrix}^2 x^n , h(x) = 
\frac{1}{1-3x} \hypergeom{2}{1}{\frac13 , \frac23}{1}{\frac{27x^2(1-x)}{(1-3x)^3}},$$
the stated equality now comes down to proving that $g(x)=h(x)$, $|x|<\frac19$.
We will show that both $g(x)$ and $h(x)$ satisfy the Heun differential equation (see \cite{Heun}) with initial values
\begin{equation}\label{heun} x(1-x)(1-9x) y'' + (1-20x+27x^2) y' + (9x-3) y =0, y(0) = 1, y'(0) = 3 . \end{equation}

If we write $g(x) = \sum_{n=0}^\infty g_n x^n$, $|x|<\frac19$, then it follows from \cite[Theorem 1; see also Table 1]{V} that
\begin{equation}\label{rec} n^2 g_n -(10n^2-10n+3)g_{n-1} + 9 (n-1)^2 g_{n-2} = 0, n\ge 2, g_0=1, g_1=3. \end{equation}
But then it is a straightforward computation that $g(x)$ satisfies \eqref{heun}. Indeed, plugging $y=g(x) = \sum_{n=0}^\infty g_n x^n$ in the left hand side of \eqref{heun} and extracting the coefficient of $x^{n-1}$ we obtain
$$ n(n-1) g_n -10(n-1)(n-2)g_{n-1} + 9(n-2)(n-3)g_{n-2} + ng_n -20(n-1)g_{n-1} +27(n-2)g_{n-2} +9g_{n-2} -3g_{n-1} = $$
$$ = n^2 g_n +(-10n^2+30n-20-20n+20-3)g_{n-1} +(9n^2-45n+54+27n-54+9)g_{n-2} = $$
$$ n^2g_n -(10n^2 -10n +3) g_{n-1} + 9(n-1)^2g_{n-2} = 0 , $$
where in the last step we use \eqref{rec}.

Next, let us turn to $h(x)$. Introduce $z(x) = \frac{27x^2(1-x)}{(1-3x)^3}$ and $w(z)= \hypergeom{2}{1}{\frac13 , \frac23}{1}{z}$. Then (see, for instance, \cite[Section 9.15]{GR}) 
$$ (1-z)zw''(z)+(1-2z)w'(z)-\frac29 w(z)= 0. $$
We have that $h(x)=\frac{1}{1-3x}w(z(x))$, 
$ h'(x) = \frac{3}{(1-3x)^2} w(z(x)) + \frac{54x}{(1-3x)^5} w'(z(x)) $, and
$$ h''(x)=\frac{18}{(1-3x)^3} w(z(x)) + \frac{54(15x+1)}{(1-3x)^6} w'(z(x)) + \frac{4(27x)^2}{(1-3x)^9} w''(z(x)) . $$  
Plugging $y=h(x)$ in the left hand side of \eqref{heun} yields
$$ x(1-x)(1-9x) h''(x) + (1-20x+27x^2) h'(x) + (9x-3) h(x)= $$
$$ \ \ \ \ \ \ \frac{108x}{(1-3x)^3} \left( (1-z(x))z(x)w''(z(x))+(1-2z(x))w'(z(x))-\frac29 w(z(x))\right)= 0. $$
In addition, it is easy to check that $h(0)=1, h'(0)=3$. 

Thus both $g(x)$ and $h(x)$ satisfy \eqref{heun}, and thus by uniqueness we find that $h(x)=g(x)$. \hfill $\square$

\begin{remark}\label{rem35}\rm
Using the Birkhoff-Trjitzinsky method (see \cite{BT}, and \cite{Immink} for complete proofs; see also \cite{WZ} and \cite{melczer}) one can obtain that the asymptotics of $g_n=h_n$ is $0.41349667\cdot \frac{9^n}{n} (1+O(n^{-1}))$. From this one can deduce that $g(x)$ is transcendental over ${\mathbb Q}(x)$; see \cite[Corollary 2.1]{melczer}.  This implies that the autoregressive filter problem in three and more variables is significantly more involved from the case of one or two variables in the sense that one can no longer expect necessary and sufficient conditions via polynomial expressions with rational coefficients, such as the low rank requirement in two variables.
\end{remark}

\begin{cor}\label{cor34}
  $$ \frac{1}{r^2-3}  \hypergeom{2}{1}{\frac13,\frac23}{1}{\frac{27(r^2-1)}{(r^2-3)^3}} = 
\frac{1}{(r-1)^{\frac32}(r+3)^{\frac12}}  \hypergeom{2}{1}{\frac12,\frac12}{1}{\frac{16r}{(r-1)^3(r+3)}} . $$
\end{cor}

{\it Proof.} Combine Theorems \ref{thmc000} and \ref{c000formula}. \hfill $\square$

There are formulas that relate $\hypergeom{2}{1}{\frac13,\frac23}{1}{z}$ and 
$\hypergeom{2}{1}{\frac12,\frac12}{1}{\widetilde{z}}$ (see, for instance, \cite[page 112]{AB}), but the above equality seems to be of a different nature than those already known.

\bigskip

We end this section by providing a proof for Theorem \ref{3varcase}.

\medskip

{\it Proof of Theorem \ref{3varcase}.} By Theorem \ref{main} we see that $\frac{|b|}{a}<1$ is necessary and sufficient. 
The proof is the same as the proof of Theorem \ref{main}, except that we will use the expression of $c_{000}$ from Theorem \ref{c000formula}. Let $d=3$. As before, the polynomial $p(z)$ after normalization so that $p(0)=1$ will satisfy \eqref{Aeq}. Starting with $A$ as in \eqref{Amatrix2} we can, by Lemma \ref{invlemma}, rescale the matrix as $\frac{a+(d-1)c}{a^2+(d-1)ac-d|b|^2} A$ so that the (1,1) entry of its inverse is 1, which corresponds to the situation where $p(0)=1$. Then, again using Lemma \ref{invlemma}, we find that $\frac{1}{r}=-\frac{dp}{dz_1}|_{z=0} $ corresponds to the value $\frac{1}{r}=\frac{b}{a^2+(d-1)ac-d|b|^2}$. Using this value for $r$ as well as
$c_{000}= \widehat{\frac{1}{|p|^2}} (0,\ldots, 0 ) = a\frac{a+(d-1)c}{a^2+(d-1)ac-d|b|^2}$, we find that Proposition \ref{prop27} yields equality \eqref{hypeq}. \hfill $\square$

\section{The three variable case: other Fourier coefficients}\label{sec4}
In \cite{GWW} the current authors considered the two variable analog, and obtained the following expression for the Fourier coefficients of $f(z_1,z_2)= |1 - \frac{z_1 +z_2}{r}|^{-2}$, $r > 2$.
\begin{theorem}\label{coeffi}  \cite[Theorem 1]{GWW}
Let $ p(z_1,z_2) = 1 - \frac{z_1 +z_2}{r} $
with $r>2$, and let $ c_{k_1,k_2}$ denote the Fourier coefficients of its spectral density function $f(z_1,z_2)= |1 - \frac{z_1 +z_2}{r}|^{-2}$.
Then we have
$$ c_{k_1,k_2} = \frac{1}{\sqrt{1-\frac{4}{r^2}}} \left(\frac{r}{2}-\sqrt{\frac{r^2}{4}-1}\right)^{|k_1|+|k_2|}, \  k_1 k_2 \le 0, $$
and 
$$ c_{k_1,k_2} = \frac{{|k_1|+|k_2| \choose |k_1|}}{r^{|k_1|+|k_2|}}\ \hypergeom{3}{2}{1,\frac{|k_1|+|k_2|}{2}+1,\frac{|k_1|+|k_2|+1}{2}}{|k_1|+1,|k_2|+1}{\frac{4}{r^2}}, \ k_1 k_2 >0. $$
\end{theorem}

In an attempt to obtain a three variable generalization of the above result, we have found following expressions for the Fourier coefficients $c_J$, $J\in \{ -1, 0 , 1\}^3$ of $f(z_1,z_2)= |1 - \frac{z_1 +z_2+z_3}{r}|^{-2}$, $r > 3$. 
\bigskip

\begin{theorem}\label{thmc000b} 
Using the same notation as in Theorem \ref{thmc000},
we have
%
%
$$c_{100}= \frac{r^2}{4\pi} \int_0^{2\pi} \frac{1}{r-e^{it}} \left( \sqrt{\frac{r^2-2r\cos t +1}{r^2-2r\cos t -3}}-1 \right) dt = $$ 
$$ -\frac{r}{2}+\frac{r^2}{4\pi}  \int_0^{2\pi} \frac{r-\cos t}{ \sqrt{r^2-2r\cos t +1}\sqrt{r^2-2r\cos t -3}} dt , $$ 
$$c_{-1,1,0}= \frac{r^2}{8\pi} \int_0^{2\pi} \sqrt{\frac{r^2-2r\cos t -3}{r^2-2r\cos t +1}} - 2 + \sqrt{\frac{r^2-2r\cos t +1}{r^2-2r\cos t -3}} dt , $$ and
 
 $$ c_{011}= \frac{r^2}{\pi}  \int_0^{2\pi} \frac{1}{e^{it}(r-e^{it})(\sqrt{(r^2+1-2r\cos t)(r^2-3-2r\cos t)}+r^2-3-2r\cos t)} dt = $$ 
$$-\frac12+ \frac{r^2}{4\pi} \int_0^{2\pi} \frac{r\cos t -\cos 2t}{\sqrt{r^2+1-2r\cos t }\sqrt{ r^2-3-2r\cos t}} dt= $$
$$ \frac{r}{4\pi} \int_0^{2\pi} \cos t \sqrt{\frac{r^2-2r\cos t +1}{r^2-2r\cos t -3}} dt +  \frac{r}{4\pi} \int_0^{2\pi} \frac{r-\cos t}{ \sqrt{r^2-2r\cos t +1}\sqrt{r^2-2r\cos t -3}} dt-\frac12 = $$
\begin{equation}\label{c110} \frac{r}{4\pi} \int_0^{2\pi} \cos t \sqrt{\frac{r^2-2r\cos t +1}{r^2-2r\cos t -3}} dt  +\frac{c_{100}}{r} . \end{equation}

 \end{theorem}

{\it Proof of Theorem \ref{thmc000b}.}
From Proposition \ref{KWformula}
we get $$c_{01} (e^{it} ) = \frac {2{r}^{2}}{r-{{\rm e}^{it}}} \left( \sqrt{(r^2+1-2r\cos t)(r^2-3-2r\cos t)}+r^2-3-2r\cos t \right) ^{-1}
$$
Using $c_{011} =\frac{1}{2\pi}  \int_0^{2\pi} c_{01}(e^{it}) e^{-it} dt $ we consequently obtain 
$$ c_{011} = \frac{r^2}{\pi} \int_0^{2\pi} \frac{1}{e^{it}(e^{it}-r)\sqrt{r^2-3-2r\cos t }(\sqrt{r^2+1-2r\cos t }+ \sqrt{ r^2-3-2r\cos t})} dt .$$
Multiplying numerator and denominator in the integrand with $\sqrt{r^2+1-2r\cos t }- \sqrt{ r^2-3-2r\cos t}$, we obtain
 $$ \frac{r^2}{4\pi} \int_0^{2\pi} \frac{\sqrt{r^2+1-2r\cos t }}{e^{it}(r-e^{it})\sqrt{ r^2-3-2r\cos t}} dt - \frac{r^2}{4\pi} \int_0^{2\pi} \frac{1}{e^{it}(r-e^{it})} dt. $$
The second term equals $\frac12$, and for the first term we can take its real part (since we know that $c_{011}$ is real).
This gives
$$ c_{011} =-\frac12+ \frac{r^2}{4\pi} \int_0^{2\pi} \frac{r\cos t -\cos 2t}{r^2+1-2r\cos t } \ \frac{\sqrt{r^2+1-2r\cos t }}{\sqrt{ r^2-3-2r\cos t}} dt = $$ $$ -\frac12+ \frac{r^2}{4\pi} \int_0^{2\pi} \frac{r\cos t -\cos 2t}{\sqrt{r^2+1-2r\cos t }\sqrt{ r^2-3-2r\cos t}} dt .$$
The last equality for $c_{011}$ is obtained by using $\frac{1}{z(r-z)}=\frac{1}{r} (\frac{1}{z} + \frac{1}{r-z})$ and applying it to the first expression for $c_{011}$.

Next, from Proposition \ref{KWformula} we find 
$$ c_{-1,1}(e^{it})=\frac{4r^2}{\sqrt{r^2+1-2r\cos t}\sqrt{r^2-3-2r\cos t}(\sqrt{r^2+1-2r\cos t} + \sqrt{r^2-3-2r\cos t})^2}. $$
Multiplying numerator and denominator with $(\sqrt{r^2+1-2r\cos t }- \sqrt{ r^2-3-2r\cos t})^2$ we obtain
$$c_{-1,1}(e^{it})= \frac{r^2}{4} \left( \frac{(\sqrt{r^2+1-2r\cos t }- \sqrt{ r^2-3-2r\cos t})^2}{\sqrt{r^2+1-2r\cos t}\sqrt{r^2-3-2r\cos t}} \right)= $$ $$
\frac{r^2}{4} \left( \sqrt{\frac{r^2-2r\cos t -3}{r^2-2r\cos t +1}} - 2 + \sqrt{\frac{r^2-2r\cos t +1}{r^2-2r\cos t -3}} \right) .$$
Use now $c_{-1,1,0} =\frac{1}{2\pi}  \int_0^{2\pi} c_{-1,1}(e^{it}) dt $  to obtain the result. 

The proof for $c_{100}$ is similar.
\hfill $\square$

\bigskip

 In Theorem \ref{thmc000} we have expressed $c_{000}$ in terms of the complete elliptic integral of the first kind.
 We can express the other Fourier coefficients above in terms of the complete elliptic integral of the first, second and third kind, which are $K(m)$, $E(m)$ and $ \Pi (n,m)$, respectively, where
$$ E(m) = \int_0^{\frac{\pi}{2}} \sqrt{1-m\sin^2 t } \ dt = \frac{\pi}{2} \hypergeom{2}{1}{-\frac12,\frac12}{1}{m}, $$ and  $$ \Pi (n,m) = \int_0^{\frac{\pi}{2}} \frac{1}{(1-n\sin^2 (t))\sqrt{1-m\sin^2(t)}} dt. $$

\begin{prop}\label{prop32}
Using the same notation as in Theorem \ref{thmc000},
we have
$$ c_{100}= \frac{r}{3}(c_{000}-1), $$
$$ c_{011}=   \frac{1}{3}(c_{000}-1) + $$ $$ \frac{(r^4-2r^2-15)K(\frac{16r}{(r+3)(r-1)^3})-(r+3)(r-1)^3E(\frac{16r}{(r+3)(r-1)^3}) -4(r-3)(r+1)\Pi (\frac{4r}{(r+3)(r-1)},
\frac{16r}{(r+3)(r-1)^3})}{4\pi(r-1)\sqrt{(r+3)(r-1)}} $$
\begin{equation}\label{c011b}\ \ \ \ \ \  = \frac{1}{3}(c_{000}-1) -\frac12  + \end{equation} $$ \frac{ (r+3)(r-1)^3K(\frac{16r}{(r+3)(r-1)^3}) -(r+3)(r-1)^3E(\frac{16r}{(r+3)(r-1)^3}) +4(r-3)(r+1)\Pi (\frac{4}{(r-1)^2}, 
\frac{16r}{(r+3)(r-1)^3})}{4\pi(r-1)\sqrt{(r+3)(r-1)}} , $$
$$ c_{111}=\frac{3}{r} c_{011}, $$
$$c_{0,1,-1}=\frac12 (rc_{001}-c_{000}), $$
$$ c_{1,1,-1}=rc_{011}-2c_{001} . $$
Other Fourier coefficients $c_J$, $J\in \{ -1, 0 , 1\}^3$, are obtained via $c_J = c_{\sigma(J)} =c_{-J}$, where $\sigma$ is a permutation. 
\end{prop}
 
 {\it Proof.}
%
First observe that 
\begin{equation}\label{fp=1/p} \frac{1}{|p(z_1, \ldots , z_d )|^2} p(z_1, \ldots , z_d) = \frac{1}{ p(\frac{1}{z_1}, \ldots , \frac{1}{z_d} )}= \sum_{k\in{\mathbb N}_0^d} \phi_k z^{-k} , z\in{\mathbb T}^d, \end{equation}
where $\phi_0=1$. If we extract the Fourier coefficients indexed by $\Lambda = \{ 0,1 \}^3$ on both sides, we obtain
$$ \begin{bmatrix} c_{000} & c_{00,-1} & c_{0,-1,0} & c_{0,-1,-1} & c_{-100} & c_{-1,0,-1} & c_{-1,-1,0} & c_{-1,-1,-1} 
\cr 
c_{001} & c_{000} & c_{0,-1,1} & c_{0,-1,0} & c_{-101} & c_{-1,0,0} & c_{-1,-1,1} & c_{-1,-1,0} 
\cr 
c_{010} & c_{01,-1} & c_{000} & c_{0,0,-1} & c_{-110} & c_{-1,1,-1} & c_{-1,0,0} & c_{-1,0,-1} 
\cr 
c_{011} & c_{010} & c_{001} & c_{000} & c_{-111} & c_{-1,1,0} & c_{-1,0,1} & c_{-1,0,0} 
\cr 
c_{100} & c_{1,0,-1} & c_{1,-1,0} & c_{1,-1,-1} & c_{000} & c_{0,0,-1} & c_{0,-1,0} & c_{0,-1,-1} 
\cr 
c_{101} & c_{100} & c_{1,-1,1} & c_{1,-1,0} & c_{001} & c_{000} & c_{0,-1,1} & c_{0,-1,0} 
\cr 
c_{110} & c_{1,1,-1} & c_{1,0,0} & c_{1,0,-1} & c_{010} & c_{0,1,-1} & c_{000} & c_{0,0,-1} 
\cr 
c_{111} & c_{110} & c_{101} & c_{100} & c_{011} & c_{010} & c_{001} & c_{000} 
\end{bmatrix} \begin{bmatrix} 1 \cr -\frac{1}{r} \cr -\frac{1}{r} \cr -\frac{1}{r} \cr 0 \cr \vdots \cr 0 \end{bmatrix} = \begin{bmatrix} 1 \cr 0 \cr \vdots \cr 0 \end{bmatrix} .$$
Since $p$ is a symmetric polynomial with real coefficients we have that $c_J = c_{\sigma(J)} =c_{-J}$, where $\sigma$ is a permutation. Thus we obtain
$$ c_{000}-\frac{3c_{001}}{r} = 1, (1-\frac{1}{r})c_{001}-\frac{c_{000}}{r}-\frac{c_{0,-1,1}}{r} =0 , c_{100}-\frac{2c_{0,-1,1}}{r} -\frac{c_{1,-1,-1}}{r} = 0, $$
$$c_{011}-\frac{c_{1,1,-1}+c_{100}+c_{0,-1,1}}{r}=0 , c_{111}-\frac{2c_{011}}{r}-\frac{c_{100}}{r} = 0 . $$
This yields the stated relations between the different Fourier coefficients (see also \cite[Proposition 3.1.1]{Wong}).

Finally, we turn to $c_{011}$. To prove the first expression for $c_{011}$, by \eqref{c110} it suffices to prove
\begin{equation}\label{intsq2}  \int_0^{2\pi} \cos t \sqrt{\frac{r^2-2r\cos t +1}{r^2-2r\cos t -3}} dt \ = \end{equation}
$$ \frac{-(r+3)(r-1)^3E(\frac{16r}{(r+3)(r-1)^3}) + (r^4-2r^2-15)K(\frac{16r}{(r+3)(r-1)^3}) -4(r-3)(r+1)\Pi (\frac{4r}{(r+3)(r-1)}, 
\frac{16r}{(r+3)(r-1)^3})}{r(r-1)\sqrt{(r+3)(r-1)}}.  $$ The left hand side of \eqref{intsq2} can be rewritten as
\begin{equation}\label{intsq5} 
\int_0^{2\pi} \cos t \sqrt{\frac{r^2-2r\cos t +1}{r^2-2r\cos t -3}} dt = \int_0^{2\pi} \sqrt{\frac{r^2-2r\cos t +1}{r^2-2r\cos t -3}} dt - \int_0^{2\pi} (1-\cos t) \sqrt{\frac{r^2-2r\cos t +1}{r^2-2r\cos t -3}} dt . \end{equation}

We will first show that 
\begin{equation}\label{intsq}  \int_0^{2\pi} \sqrt{\frac{r^2-2r\cos t +1}{r^2-2r\cos t -3}} dt =\frac{4\left( 4K(\frac{16r}{(r+3)(r-1)^3})+(r-3)(r+1) \Pi(\frac{4r}{(r+3)(r-1)} , \frac{16r}{(r+3)(r-1)^3})\right)}{(r-1)\sqrt{(r+3)(r-1)}} .  \end{equation}
To prove \eqref{intsq} we need to show the equality
$$  \int_0^{2{\pi}} \sqrt{\frac{r^2-2r\cos t +1}{r^2-2r\cos t -3}} dt =
 \int_0^{2{\pi}} \frac{\sqrt{(r+3)(r-1)^3-16r\sin^2 t}}{(r+3)(r-1)-4r\sin^2 t}dt
 .
$$
To prove the above equality we make some simplifications. In the second integral because everything is in terms of $\sin^2 t$ the integral from $[0,2\pi]$ is equal to 4 times the integral $[0,\pi/2]$. For the first integral make the change of variables $\cos t=1-2\sin^2 t/2$ then $t\to t/2$ then put everything on $[0,\pi/2]$ and divide by 4 to obtain,
$$
\frac{r-1}{\sqrt{(r-3)(r+1)}}\int_0^{\frac{\pi}{2}} \sqrt{\frac{1+\frac{4r}{(r-1)^2}\sin^2 t}{1+\frac{4r}{(r-3)(r+1)}\sin^2 t}} dt =
\sqrt{\frac{r-1}{r+3}} \int_0^{\frac{\pi}{2}} \frac{\sqrt{1-\frac{16r}{(r+3)(r-1)^3}\sin^2 t}}{1-\frac{4r}{(r+3)(r-1)}\sin^2 t}dt.
$$
Now let $p^2=-\frac{4r}{(r-1)^2}$ and $q^2=-\frac{4r}{(r-3)(r+1)}$. Then $1-q^2=\frac{(r+3)(r-1)}{(r-3)(r+1)}$ and 
$\frac{q^2-p^2}{1-q^2}=\frac{16r}{(r+3)(r-1)^3}$. The integrals become
\begin{equation}\label{intsq6} 
\int_0^{\frac{\pi}{2}} \sqrt{\frac{1-p^2\sin^2 t}{1-q^2\sin^2 t}} dt =
\frac{1}{\sqrt{1-q^2}} \int_0^{\frac{\pi}{2}} \frac{\sqrt{1-\frac{q^2-p^2}{1-q^2}\sin^2 t}}{1+\frac{q^2}{1-q^2}\sin^2 t}dt.
\end{equation}
On the right hand integral make the change of variable $\sin t=\frac{\sqrt{1-q^2}\sin x}{\sqrt{1-q^2\sin^2x}}$\footnote{The reverse change of variables is $\sin x = \frac{\sin t}{\sqrt{1-q^2+q^2\sin^2 t}}$ and we get
$ dx =\frac{\sqrt{1-q^2}}{1-q^2+q^2\sin^2 t}dt. $}, then $[0,\pi/2]\to[0,\pi/2]$ and the right hand integral goes to the left hand integral.
Indeed, we have $$ \cos t dt = \frac{\sqrt{1-q^2\sin^2 x}\sqrt{1-q^2}\cos x - \sqrt{1-q^2}\sin x \frac{1}{2\sqrt{1-q^2\sin^2 x}}(-2q^2\sin x \cos x)}{1-q^2\sin^2 x} dx . $$
Using $\cos t = \sqrt{1 - \frac{(1-q^2)\sin^2x}{1-q^2\sin^2 x}} = \frac{\cos x}{\sqrt{1-q^2\sin^2 x}}$, we find
$$ dt = \frac{\sqrt{1-q^2}}{1-q^2\sin^2 x}dx. $$
Now equality \eqref{intsq6} (and thus \eqref{intsq}) follows after some manipulations.

Next, we deal with the second term of the right hand side of \eqref{intsq5}:
\begin{equation}\label{intsq7} \int_0^{2\pi} (1-\cos t) \sqrt{\frac{r^2-2r\cos t +1}{r^2-2r\cos t -3}} dt = 2
\int_0^{\pi}(1-\cos t) \sqrt{\frac{\frac{r^2-1}{2r}-\cos t }{\frac{r^2-3}{2r}-\cos t }} dt. \end{equation}
By using the change of variables $u= \cos t$ (and thus $dt = -\frac{1}{\sqrt{1-u^2}} du$), we can rewrite this as
$$ 2 \int_{-1}^1 \sqrt{\frac{(\frac{r^2+1}{2r} - u)(1-u)}{(\frac{r^2-3}{2r} - u)(u-(-1))}} du . $$ 
Let $$ a=\frac{r^2+1}{2r} , b=\frac{r^2-3}{2r}, c=1, y=1, d=-1 , $$ and observe that
$a>b>c\ge y>d$. 
We can now use \cite[Equations 252.17 and 362.16]{Byrd}, which yield
$$ \int_d^y \sqrt{\frac{(c-u)(a-u)}{(b-u)(u-d)}} du =   $$
$$ \frac{(a-d)(c-d)g}{2\alpha^2(k^2-\alpha^2)} \left( \alpha^2E(k^2)+(k^2-\alpha^2)K(k^2)+(2k^2\alpha^2-\alpha^4-k^2)\Pi(\alpha^2,k^2)\right), $$
where 
$$g=\frac{2}{\sqrt{(a-c)(b-d)}}, 
\alpha^2=\frac{d-c}{a-c},
k^2=\frac{(a-b)(c-d)}{(a-c)(b-d)}. $$
We obtain that \eqref{intsq7} equals
\begin{equation}\label{intsq8} \frac{(r-1)^3(r+3)E(\frac{16r}{(r+3)(r-1)^3})-(r-1)^2(r+1)^2K(\frac{16r}{(r+3)(r-1)^3})+4(r+1)^3\Pi(\frac{-4r}{(r-1)^2},\frac{16r}{(r+3)(r-1)^3})}{r \sqrt{(r-1)^3(r+3)}}. \end{equation}
Next we observe that \cite[Equation 117.03]{Byrd}, after multiplying with $\frac{(r-1)^5(r+3)}{4r(r+1)^2}$, gives
\begin{equation}\label{intsq9} (r+1)^2 \Pi (\frac{-4r}{(r-1)^2},\frac{16r}{(r+3)(r-1)^3}) = \end{equation} $$(r-3)(r+1) \Pi (\frac{4r}{(r+3)(r-1)}, 
\frac{16r}{(r+3)(r-1)^3}) +4 K (\frac{16r}{(r+3)(r-1)^3}) . $$
Putting these together with \eqref{intsq}, yields \eqref{intsq2}.


To prove the second equality for $c_{011}$ from the first, we use (see \cite[Formula 117.02]{Byrd})
$$ \Pi(n,m)= K(m)-\Pi(\frac{m}{n}, m) + \frac{\pi}{2}\sqrt{\frac{n}{(1-n)(n-m)}} , $$
with $n=\frac{4r}{(r+3)(r-1)}$ and $m=\frac{16r}{(r+3)(r-1)^3}$.  
The constant here works out to equal
$\frac{\pi}{2} \frac{(r+3)^{\frac12}(r-1)^{\frac32}}{(r-3)(r+1)}. $
Thus \eqref{c011b}
follows. \hfill $\square$

\medskip
Equation \eqref{fp=1/p} yields the relations
$$ c_{klm} - \frac{c_{k-1,l,m}+c_{k,l-1,m}+c_{k,l,m-1}}{r} =0 , (k,l,m)\not\in -{\mathbb N}_0^3. $$
These equalities provide a partial picture of the Fourier coefficients of  $|1 - \frac{z_1 +z_2+z_3}{r}|^{-2}$, $r>3$. 
Our method to determine other relations rely on the formulas obtained in Proposition \ref{KWformula}. The inverses in this proposition are obtained via \cite[Theorem 1.1]{KW} and the ability to find a formula for the inverse of a tridiagonal infinite Toeplitz matrix. If we want to use this method to obtain expressions for Fourier coefficients beyond the indices $\{-1,0,1\}^3$, we will need to be able to find manageable expressions for (part of) the inverse of more involved infinite (block) Toeplitz matrices, which is a challenge.

\bigskip
{\bf Acknowledgment.} We would like to thank Professor Bruce C. Berndt for his feedback on Corollary \ref{cor34}. Also we thank Paul D. Hanna for giving us background on his entries on the On-Line Encyclopedia of Integer Sequences. Finally, we thank Professor Stephen Melczer for his help with Remark \ref{rem35}.

\bibliographystyle{plain}

\begin{thebibliography}{10}
%
\bibitem{AB} George E. Andrews and Bruce C. Berndt, {\em Ramanujan's lost notebook. Part V.} Springer Verlag, Cham. xii, 430 p. (2018). 

\bibitem{BT} George D. Birkhoff and W. J. Trjitzinsky. Analytic theory of singular difference equations. {\em Acta Math.}, 60(1):1--89, 1933.

\bibitem{Boch}
S.~Bochner.
\newblock Vorlesungen \"{u}ber {F}ourierische {I}ntegral.
\newblock Akademische Verlagsgesellschaft, Leipzig, 1932.

\bibitem{Boch1}
S.~Bochner.
\newblock Monotone {F}unktionen, {S}tieltjesschee {I}ntegrale und harmonische
  {A}nalyse.
\newblock {\em Math. Anal.}, 108:378--410, 1933.

\bibitem{BW}
Charles Burnette and Chung Y. Wong, {Abelian Squares and Their Progenies}, preprint.

\bibitem{Byrd} Paul F. Byrd and Morris D. Friedman, {\em Handbook of Elliptic Integrals for Engineers and Scientists}, Springer-Verlag, Berlin Heidelberg, 1971.

\bibitem{GW}
Jeffrey~S. Geronimo and Hugo~J. Woerdeman.
\newblock Positive extensions, {F}ej\'er-{R}iesz factorization and
  autoregressive filters in two variables.
\newblock {\em Ann. of Math. (2)}, 160(3):839--906, 2004.


\bibitem{GWW} Jeffrey S. Geronimo, Hugo J. Woerdeman and Chung Y. Wong, Spectral density functions of bivariable stable polynomials, preprint. arXiv:2012.12980.

\bibitem{GR} Izrail Solomonovich Gradshteyn, Iosif Moiseevich Ryzhik, Yuri Veniaminovich Geronimus, Michail Yulyevich Tseytlin, Alan Jeffrey, Daniel Zwillinger, Victor Hugo Moll (eds.). {\em Table of Integrals, Series, and Products}. Translated by Scripta Technica, Inc. (8 ed.), 2015. Academic Press, Inc. ISBN 978-0-12-384933-5

\bibitem{Hanna} Paul D. Hanna, contribution dated February 26, 2012, regarding the  sequence A002893 on the On-Line Encyclopedia of Integer Sequences (oeis.org).

\bibitem{Hanna2} Paul D. Hanna, private communication.

\bibitem{HL}
Henry Helson and David Lowdenslager, Prediction theory and
{F}ourier
  series in several variables, {\em Acta Math.} {99} (1958), 165--202. 

\bibitem{HL2}
Henry Helson and David Lowdenslager, Prediction theory and {F}ourier series in several variables. II, {\em Acta Math.} 106 (1961), 175--213.

\bibitem{Heun} K. Heun. Zur Theorie der Riemann'schen Functionen zweiter Ordnung mit vier Verzweigungspunkten. {\em Math. Ann.}, 33(2):161--179, 1888. 

\bibitem{Immink} G. K. Immink, 
Reduction to Canonical Forms and the Stokes Phenomenon in the Theory of Linear Difference Equations
{\em SIAM J. Math. Anal.}, 22(1), 238--259, 1991.

\bibitem{Kailath1986}
Thomas Kailath.
\newblock A theorem of {I}. {S}chur and its impact on modern signal processing.
\newblock In {\em I. {S}chur methods in operator theory and signal processing},
  volume~18 of {\em Oper. Theory Adv. Appl.}, pages 9--30. Birkh\"auser, Basel,
  1986.

\bibitem{Kayran1996}
Ahmet~H. Kayran.
\newblock Two-dimensional orthogonal lattice structures for autoregressive
  modeling of random fields.
\newblock {\em IEEE Trans. Signal Process.}, 44(4):963--978, 1996.


\bibitem{KW} Selcuk Koyuncu and Hugo J. Woerdeman,  The Inverse of a Two-level Positive Definite Toeplitz Operator Matrix. In: H. Dym, M. A. Kaashoek, P. Lancaster, H. Langer, L. Lerer (eds), A Panorama of Modern Operator Theory and Related Topics. {\em Operator Theory: Advances and Applications}, vol 218 (2012), Springer, Basel.

\bibitem{Lev-Ari1989}
Hanoch Lev-Ari, Sydney~R. Parker, and Thomas Kailath.
\newblock Multidimensional maximum-entropy covariance extension.
\newblock {\em IEEE Trans. Inform. Theory}, 35(3):497--508, 1989.

\bibitem{melczer} Stephen Melczer, {\em An Invitation to Analytic Combinatorics: From One to Several Variables}, Springer Texts and Monographs in Symbolic Computation, 2020.

\bibitem{RR} L. B. Richmond and C. Rousseau, Comment on problem 87-2, SIAM Review, 31, (1989),
no. 1, 122--€"125.

\bibitem{RS} L. B. Richmond and Jeffrey Shallit, Counting Abelian Squares, The Electronic Journal of Combinatorics 16(1), 2009, R72.

\bibitem{V} H. A. Verrill, Sums of squares of binomial coefficients, with applications to Picard-Fuchs equations, ArXiv:math/0407327

\bibitem{WZ} Jet Wimp and Doron Zeilberger. Resurrecting the asymptotics of linear recurrences. {\em J. Math.
Anal. Appl.}, 111(1):162--176, 1985.

\bibitem{Wong} Chung Y. Wong, Spectral Density Functions and Their Applications. Thesis (Ph.D.)--Drexel University. 2016. 



\end{thebibliography}
 
\end{document}